\documentclass[runningheads]{llncs}

\usepackage[table]{xcolor}
\definecolor{lightgray}{gray}{0.9}
\usepackage{enumitem}
\usepackage{amssymb, amsmath}
\usepackage{amsrefs}
\usepackage{hyperref}
\usepackage{cleveref}

\usepackage{cases}

\newcommand{\Z}{\mathbb{Z}}
\newcommand{\N}{\mathbb{N}}

\newcommand{\ab}{\textsf{ab}}

\newcommand{\mbb}{\mathbb}

\usepackage{tikz}
\usepackage[all]{xy} 
\usepackage{bm}

\newcommand{\bzero}{\bm{0}}

\newcommand{\Parikh}[1]{\Psi(#1)}
\newtheorem{thmletter}{Theorem}

\begin{document}
\title{Word equations, constraints, and formal languages}

\author{Laura Ciobanu 
%\inst{Corresponding author}
\orcidID{0000-0002-9451-1471}} %
\authorrunning{L. Ciobanu}

\institute{Heriot-Watt University and Maxwell Institute,  Edinburgh EH14 4AS,
 Scotland\\
\email{L.Ciobanu@hw.ac.uk}}

\maketitle

\begin{abstract}
In this short survey we describe recent advances on word equations with non-rational constraints in groups and monoids, highlighting the important role that formal languages play in this area.


\end{abstract}

\section{Introduction}
  The question of finding solutions to equations in algebraic structures such as semigroups, groups or rings is a fundamental topic in mathematics and theoretical computer science that finds itself at the frontier between decidability and undecidability. Beyond asking about the existence of solutions to an equation it is often of interest, for both practical and theoretical reasons, to ask that the solutions  belong to a certain set or satisfy a certain property; in that case one asks about the decidability of solving \emph{equations with constraints}. 
  
  In the context of monoids and groups there are broadly two types of constraints: \emph{rational} and \emph{non-rational}. The rational ones concern properties of the solutions that can be recognised by finite automata, e.g. if there are solutions of even length in a given monoid or group, with respect to a fixed alphabet (or generating set). These are well understood in most settings, and typically in the cases when the satisfiability of equations is decidable, adding rational constraints or versions thereof remains decidable (\cite{dgh01, DG}). 
  
  The second class of constraints, the \emph{non-rational} ones, refers to sets that cannot be recognised by finite automata; examples of such constraints are relations between the lengths of the solutions or their Parikh images, or membership in context-free sets.  There is an extensive literature on word equations, that is, equations in free monoids, with various non-rational constraints \cite{Abdulla, RichardBuchi1988, DayManeaWE, GarretaGray, majumdar2021quadratic} %{day2018satisfiability} 
 and, recently, equations in groups with similar (non-rational) constraints have received attention. We highlight here some of the new results in both free monoids and groups, and refer the reader to \cite{CiobanuGarreta, CiobanuEvettsLevine, CiobanuLevine, CiobanuZetzsche2024a} for a more in-depth treatment.

 \section{Word equations with constraints}\label{sec:wordeqs}
 
 In this section we discuss equations in free monoids, which are classically referred to as `word equations' in the literature and are defined as follows.
 
 Let $\Omega$ be a finite collection of
variables and $A$ a finite set of letters, and consider  the word equation $U=V$ over $(A,\Omega)$, where $U, V$ are words in the free monoid $(A \cup \Omega)^*$.
% is a pair $(U,V)$ with $U,V\in (A\cup\Omega)^*$. Here, $U$ and $V$ are the two sides
%of an equation and we are looking to replace the variable occurrences in $U$
%and $V$ so that the two sides become equal.
A \emph{solution} to $U=V$  is a morphism
$\sigma\colon(\Omega \cup A)^*\to A^*$ that fixes $A$ point-wise and such that $\sigma(U)=\sigma(V)$. We will refer to the question of whether a given equation has solutions as (\textbf{WordEqn}). This problem is solvable by the groundbreaking work of Makanin \cite{mak83a} (and subsequent improvements).
 
 %\mysubsection{Word equations}

A \emph{word equation with rational constraints} is a word equation $U=V$
together with a regular language $R_X\subseteq A^*$ for each variable
$X\in\Omega$. In this situation, we say that $\sigma\colon\Omega\to A^*$ is a
\emph{solution} if it is a solution to the word equation $U=V$ and also
satisfies $\sigma(X)\in R_X$ for each $X\in\Omega$. We will refer to the question of whether a given equation has solutions satisfying rational constraints as (\textbf{WordEqn, RAT}). This is decidable by \cite{dgh01}.

To describe the relevant non-rational constraints in this paper we need to introduce a few basic concepts. Let $|w|$ denote the word length of any $w \in A^*$. For any $a\in A$, let $| w |_a$ count the number of occurrences of $a$ in $w$; for example, $|abab^2|_a=2$. Denote by $\Parikh{L} \subseteq \N^k$ (with $k=|\Sigma|$), the Parikh image of a set $L \subseteq A^*$, that is, the abelianisation of $L$.

The starting point of the study of word equations with non-rational constraints is B\"uchi and Senger's well-known paper \cite{RichardBuchi1988}, where they show that (positive) integer addition and multiplication can be encoded into word equations in a free monoid, when requiring that the solutions satisfy one of several counting constraints; by the undecidability of Hilbert's 10th problem, deciding the satisfiability of such enhanced word equations is also undecidable.

We present the problems concerning word equations with non-rational constraints via an easy example, where each question represents one of the types of constraints that we are interested in. 

\begin{example}
Let $A=\{a, b\}$, $\Omega=\{X, Y, Z\}$ and consider the word equation 
\begin{equation}\label{eqn}
X^2aY=YZ^2a.
\end{equation}
A possible solution is $\sigma_0=(ab, a, ba)$, as $\sigma_0(X)=ab$, $\sigma_0(Y)=a$ and $\sigma_0(Z)=ba$ give $\sigma_0(X^2aY)=\sigma_0(YZ^2a)$.
\begin{enumerate}
\item (\textbf{WordEqn, LEN}) asks about the existence of solutions whose lengths satisfy some given linear system of equations over the integers.

\medskip

 For example, is there a solution $\sigma$ to (\ref{eqn}) such that $$|\sigma(X)| = |\sigma(Z)|+1?$$
A simple argument comparing the lengths of solutions on the left and right hand side will show that no such solution exists.

\medskip

\item (\textbf{WordEqn, EXP-SUM}) asks about the existence of solutions where the number of occurrences/exponent-sums of letters from $A$ satisfy some given linear system of equations over the integers.

\medskip

For example, is there a solution $\sigma$ to (\ref{eqn}) such that $$|\sigma(X)|_a = 2|\sigma(Z)|_a?$$

A simple argument comparing the numbers of $a$'s on the left and right hand side will show that no such solution exists.

\medskip

\item[3.] (\textbf{WordEqn, PARIKH}) asks about the existence of solutions whose Parikh images satisfy some given linear system of equations over the integers. (These are called $\textsf{AbelianEq}$ in \cite{DayManeaWE}.)

\medskip

For example, is there a solution $\sigma$ such that $$\Psi(\sigma(X)) = \Psi(\sigma(Z))?$$

The answer is clearly `yes', as solution $\sigma_0$ shows: $\Psi(ab)=\Psi(ba)=(1,1).$

\medskip

\item[4.] (\textbf{WordEqn, INEQ}) asks about the existence of solutions for which length, exponent-sum or other counting identities never hold.

\medskip

\begin{itemize}

\item[(i)] For example, is there a solution $\sigma$ such that $$|\sigma(X)|_a \neq |\sigma(Z)|_a?$$
If we count the possible numbers of $a$'s on the left and right hand sides of the equation, and consider the fact that both $X$ and $Z$ appear with the same exponent, we get $|\sigma(X)|_a \neq |\sigma(Z)|_a$ for every solution $\sigma$; thus the answer to the question above is negative.

\medskip

\item[(ii)] Another such problem would be: Is there a solution $\sigma$ such that $$|\sigma(X)| \neq |\sigma(Z)|?$$
An immediate length argument shows that the answer is negative.

\end{itemize}
\end{enumerate}

\end{example}

As B\"{u}chi and Senger showed~\cite[Corollary 4]{RichardBuchi1988},
problems such as (2) and (3) in full generality are undecidable. In fact, questions (2) and (3) are equivalent (see \cite[Lemma 3.2]{CiobanuGarreta}). However, it remains a well-known open problem
whether the satisfiability of word equations with \emph{length constraints}, as in (1), is decidable. Deciding algorithmically whether a word equation has solutions satisfying certain linear length constraints is a major open question, and it has deep implications, both theoretical (if undecidable, it would offer a new solution to Hilbert's 10th problem about the satisfiability of polynomial equations with integer coefficients) and practical, in the context of string solvers for security analysis. We refer the reader to the surveys \cite{Amadini,ganesh} for an overview of the area from several viewpoints, of both theoretical and applied nature. 

One of the recent results about counting constraints concerns questions such as (4), where one is asking about certain properties of the solutions to \emph{never} hold. The result below applies to more general counting functions than just length and exponent-sums of letters in the solutions. 

\begin{theorem}\cite[Theorem 5]{CiobanuZetzsche2024a}\label{mainGeorg} The problem of satisfiability of word equations with rational constraints and counting inequations is decidable.
\end{theorem}

This result relies on two important ingredients, both of which revolve around formal languages:
\begin{itemize}
\item[1.] Solution sets to word equations, as represented in (\ref{eqnrep}), are EDT0L, and therefore indexed. (\cite{CiobanuDiekertElder2016ijac})

\item[2.] A certain semilinear set, called the \emph{slice closure}, is effectively computable for any indexed language $L$. Moreover, this set, containing the Parikh image of $L$, records the counting properties (such as lengths of words, exponent-sums of letters etc.) of $L$ to a large degree. (\cite[Section 2]{CiobanuZetzsche2024a})
\end{itemize}

\newcommand{\enc}[1]{\mathsf{enc}(#1)}
\newcommand{\cT}{\mathcal{T}}
%\subsection{Representing solutions} 
We represent the set of solutions to word equations as follows. Given a map $\sigma\colon\Omega\to A^*$ with $\Omega=\{X_1,\ldots,X_k\}$, we define
\begin{equation}\label{eqnrep}
\enc{\sigma}~ = ~\sigma(X_1)\#\cdots\#\sigma(X_k), 
\end{equation}
where $\#$ is a fresh letter not in $A$. Hence, $\enc{\sigma}\in(A\cup\{\#\})^*$.
Let $\Sigma=A\cup\{\#\}$. 

The counting constraints in Theorem \ref{mainGeorg} can be expressed using rational counting functions; these are functions
$f\colon \Sigma^*\to \Z^n$ for which the graph of $f$ is a rational subset of  $\Sigma^*\times\Z^n$, that is, $\{(w,f(w)) \mid w\in\Sigma^*\}\subseteq\Sigma^*\times\Z^n$ can be defined by a rational transducer. 
%A \emph{counting transducer} is a finite automaton $\cT=(Q,\Sigma,E,q_0,F)$, where $Q$ is a finite set of \emph{states}, $\Sigma$ is its \emph{input alphabet}, $E\subseteq Q\times\Sigma^*\times\Z^n\times Q$ is its \emph{edge relation}, $q_0\in Q$ is its \emph{initial state}, and $F\subseteq Q$ is its set of \emph{final states}. A \emph{run} is a sequence
%\[ q_0(w_1,\bx_1)q_1\cdots (w_m,\bx_m)q_m \]
%with $q_0,\ldots,q_m\in Q$, $w_1,\ldots,w_m\in\Sigma^*$, and $\bx_1,\ldots,\bx_m\in\Z^n$ such that each $(q_i,w_{i+1},\bx_{i+1},q_{i+1})\in E$, and $q_m\in F$. The relation \emph{defined} by the counting transducer, denoted $R(\cT)$ is the set of all pairs $(w_1\cdots w_m,\bx_1+\cdots+\bx_m)\in \Sigma^*\times\Z^n$ for runs as above. The counting transducer is called \emph{functional} if for every word $w\in\Sigma^*$, there is at most one vector $\bx\in \Z^n$ such that $(w,\bx)\in R(\cT)$. A function $f\colon \Sigma^*\to\Z^n$ is called a \emph{counting function} if its graph $\{(w,f(w)) \mid w\in\Sigma^*\}\subseteq\Sigma^*\times\Z^n$ is defined by a functional counting transducer. 
These counting functions can then be applied to encodings of solutions, as in the following examples from \cite{CiobanuZetzsche2024a}.
\begin{enumerate}
\item Exponent-sums/letter occurrence counting: 

\medskip

The function $f_{X,a}$ with $f_{X,a}(\enc{\sigma})=|\sigma(X)|_a$ for some letter $a\in A$ and variable $X\in\Omega$. Here,  the transducer increments the counter for each $a$ between the $i$-th and $(i+1)$-st occurrence of $\#$, where $X=X_i$.

\medskip

\item Counting positions with MSO properties:

\medskip

 Consider a monadic second-order logic (MSO) formula $\varphi(x)$ with one free first-order variable $x$, evaluated over finite words. Then one can define the function $f_{\varphi}$ such that $f_{\varphi}(\enc{\sigma})$ is the number of positions $x$ in $\enc{\sigma}$ where $\varphi(x)$ is satisfied. Then $f_{\varphi}$ is a counting function, which follows from the fact that MSO formulas define regular languages. For example, we could count the number of $a$'s such that there is no $c$ between the $a$ and the closest even-distance $b$.

\medskip

\item Linear combinations of counting functions:

\medskip

 If $f_1,\ldots,f_m\colon\Sigma^*\to\Z^n$ are counting functions, then so is $f$ with $f(w)=\lambda_1f_1(w)+\cdots+\lambda_mf_m(w)$ for some $\lambda_1,\ldots,\lambda_m\in\Z$. This can be shown with a simple product construction.
 
 \medskip
 
\item\label{length-function} Length functions: 

\medskip

The function $L_X$ with $L_X(\enc{\sigma})=|\sigma(X)|$ for some $X\in\Omega$. For this, we can take a linear combination of letter counting functions.
\end{enumerate}

Using rational counting functions one can express the problems (\textbf{WordEqn, LEN}), (\textbf{WordEqn, LEN}) or (\textbf{WordEqn, LEN}) in the general framework:

%The following is the problem of \emph{word equations with rational and counting constraints}:
\begin{description}
\item[Input] A word equation $(U,V)$ with rational constraints and a counting function $f\colon (A\cup\{\#\})^*\to\Z^n$
\item[Question] Is there a solution $\sigma$ that satisfies $f(\enc{\sigma})=0$?
\end{description}
%We say that $f\colon (A\cup\{\#\})^*\to\Z^n$ is a \emph{length
%constraint} if each projection to a component of $\Z^n$ is a linear combination
%of length functions as in \cref{length-function} above. Thus, a length
%constraint can express that $\sigma(X)$ has the same length as $\sigma(Y)$ for
%$X,Y\in\Omega$. The problem of \emph{word equations with rational and length
%constraints} is the following:
%\begin{description}
%\item[Given] A word equation $(U,V)$ with rational constraints and a length constraint $f\colon (A\cup\{\#\})^*\to\Z^n$
%\item[Question] Is there a solution $\sigma$ that satisfies $f(\enc{\sigma})=\bzero$?
%\end{description}
%
%We will show here that a restriction of word equations with rational
%constraints and counting constraints is decidable. In this restriction, we
%require $f(\enc{\sigma})\ne\bzero$ rather than $f(\enc{\sigma})=\bzero$. 

Then the question answered positively in Theorem \ref{mainGeorg}, that is, the
problem of solving \emph{word equations with rational constraints and counting
inequations}, can be phrased as follows.
\begin{description}
\item[Input] A word equation $(U,V)$ with rational constraints and a counting function $f\colon (A\cup\{\#\})^*\to\Z^n$
\item[Question] Is there a solution $\sigma$ that satisfies $f(\enc{\sigma})\ne \bzero$?
\end{description}
%Here, we show that word equations with rational constraints and counting inequations are decidable (\cref{counting-inequations-decidable}).

% The question of decidability of (systems of) equations is widely known as the
%{\em Diophantine Problem} for $G$, and we denote it by $\mc{DP}(G)$.

\section{Groups equations with constraints}

 We find similarities to word equations with constraints, but also new territory, when entering the world of groups. The question of decidability of (systems of) equations is widely known as the {\em Diophantine Problem} in groups rather than solving {\em Word Equations}, but to have a uniform notation throughout the paper we will denote this question by (\textbf{WordEqnGp}) rather than the Diophantine Problem.% asking about the existence of algorithms to determine  whether a given equation has solutions in a group as .
 
 The main classes of groups (not necessarily disjoint, not an exhaustive list), for which (\textbf{WordEqnGp}) is decidable, are listed below:
 
 \begin{enumerate}
 \item virtually abelian groups \cite{EvettsLevine},
 \item hyperbolic groups, with subclasses treated in several stages, chronologically:
 \begin{enumerate}
 \item free groups \cite{mak83a},
 
 \item torsion-free hyperbolic groups \cite{RS95},
 
 \item virtually free and (torsion) hyperbolic groups \cite{DG},
 \end{enumerate}
 
 \item partially commutative groups (or right-angled Artin groups)\cite{DiekertMuscholl},
 
 \item graph products of groups with decidable (\textbf{WordEqnGp}) \cite{DLijac},
 \item virtually direct products of hyperbolic groups (which include dihedral Artin groups, and therefore all Baumslag-Solitar groups $BS(m,m)$) \cite{equations_VDP},
 
 \item central extensions of hyperbolic groups \cite{liang_central_extensions}, 
 
 \item some torsion-free relatively hyperbolic groups \cite{Dahmani}. 
 \end{enumerate}
 
 For several additional classes of groups there are positive results for restricted types of equations, for example:
 
 \begin{enumerate}
 \item single equations in the Heisenberg group \cite{DuchinLiangShapiro},
 
 \item quadratic equations in the Baumslag-Solitar groups $BS(1,n)$ \cite{MandelUsh},
 
 \item quadratic equations in the Grigorchuck group \cite{quadGrig}.
 \end{enumerate}
 
 In many of the above cases the decidability of (\textbf{WordEqnGp}) with rational constraints, so (\textbf{WordEqnGp, RAT}), has also been established. Notable exceptions are the partially commutative groups and the hyperbolic groups, where (\textbf{WordEqnGp, RAT}) is undecidable in full generality, and one needs to weaken the type of constraints to get decidability.
 
 To make this paper accessible to a wide audience we give some background before describing the non-rational constraints that are considered in groups.
 
 \subsection{Background on group theory}
% \subsection{Preliminaries}
%\label{sec:Prelim}
Let $\Sigma$ be a finite set, and as before, let $|w|$ denote the word length of any $w \in \Sigma^*$. For any $a\in \Sigma$, let $| w |_a$ count the number of occurrences of $a$ in $w$.%; for example, $||abab^2||_a=2$.
%that is, $||w||_a$ is the number of occurrences of $a$ in the word $w\in \Sigma^*$. 

\textbf{Free groups.} Define $\Sigma^{-1}$ as the set of formal inverses of elements in $\Sigma$ and denote the free group with generating set $\Sigma$ by $F(\Sigma)$; $F(\Sigma)$ can be viewed as the set of all \emph{freely reduced words} over $\Sigma^{\pm1}=\Sigma \cup \Sigma^{-1}$, that is, words not containing $xx^{-1}$ as subwords, $x\in \Sigma^{\pm1}$, together with the operations of concatenation and free reduction (that is, the removal of any $xx^{-1}$ that might occur when concatenating two words, where $(x^{-1})^{-1}$ for any $x\in\Sigma^{\pm 1}$.

\textbf{Partially commutative groups (RAAGs)}. A group $G$ generated by the finite set $\Sigma$, subject to a set of relations $R \subset F(\Sigma) \times F(\Sigma)$, is denoted as $G=\langle \Sigma \mid R \rangle$ and can be viewed as $F(\Sigma)$ modulo the relations in $R$: two elements are equal in $G$ if there is a way of writing them as words which can be transformed into each other via the relations in $R$, together with deleting or inserting  $xx^{-1}$, $x\in\Sigma^{\pm}$. For example, the free abelian group $(\mbb{Z}^2, +)$ can be given by generators $\Sigma=\{a,b\}$ that satisfy the relation $(ab,ba)$, which we write as $ab=ba$. One may replace $ab=ba$ by $aba^{-1}b^{-1}=1$, and use the commutator notation $[a,b]=1$. 

A class of groups that lie between the free (non-abelian) and the free abelian groups, in terms of their presentations, are the \emph{partially commutative groups}. They are the group theoretic counterpart to partially commutative monoids, or trace monoids. In geometric group theory these are called \emph{right-angled Artin groups} (RAAGs), and we will use this short-hand notation to save space. The most common way of describing a RAAG is via a finite undirected graph $\Gamma$ with no auto-adjacent vertices (i.e. no loops at any vertex) and no multiple edges between two vertices, and letting the vertices of $\Gamma$ be the generators of the RAAG $G\Gamma$ based on $\Gamma$. The relations between generators correspond to the edges: for every edge $(u,v)$ in $\Gamma$ we introduce the commuting relation $uv=vu$.

One often writes $G$ instead of $G\Gamma$ when the graph $\Gamma$ is unambiguous.
%Example 1 below shows a RAAG based on a graph with vertices $\{a,b,c,d\}$ and $3$ edges which give the relations $ab=ba$ and $bc=cb$, and $cd=dc$. It is customary to replace $ab=ba$ by $aba^{-1}b^{-1}=1$, and write this in commutator notation, that is, $[a,b]=1$.
\begin{example}\label{ex:3edgelinesegment}
Let $\Gamma$ be the graph below with vertices $\{a,b,c,d\}$ and $3$ edges. The RAAG $G=G\Gamma$ based on $\Gamma$ has the presentation $\langle a,b,c,d \mid [a,b]=[b,c]=[c,d]=1 \rangle$; that is, $G$ has generators $\{a,b,c,d\}$ which satisfy the relations $ab=ba$, $bc=cb$, and $cd=dc$.
 \[
  \xymatrix{\stackrel{a}{\bullet}\ar@{-}[r] &\stackrel{b}{\bullet} \ar@{-}[r]&\stackrel{c}{\bullet}\ar@{-}[r] &\stackrel{d}{\bullet} }.
 \]

 \end{example}

\textbf{The abelianization of a group.} 
The abelianization of a group is the analogue of the Parikh map $\Psi: \Sigma^* \mapsto \mbb{N}^{|\Sigma|}$ for free monoids, already defined, but where instead of $\Sigma^*$ we consider a group defined by $\Sigma$.

For any group $G$, let $\ab:G\to G^{\ab}$ be the natural abelianization map. %(we will often use the additive notation $+$ in the abelian group $G/G'$). 
That is, $G^{\ab}$ is the group with the presentation of $G$, plus the additional relations that any two generators commute. In algebraic terms, $G^{\ab}= G/G'$, that is, the quotient of $G$ by its commutator subgroup; $G^{\ab}$ is a commutative group which will decompose into an infinite part of the form $\mathbb{Z}^m$, for some $m\geq 1$, and a finite abelian group $H$. The integer $m$ is called \emph{the free rank} of $G^{\ab}$. So for example, the abelianization of $F(a,b)$ is the free abelian group $(\mbb{Z}^2, +)$ and is of free rank $2$; this is because $F(a,b)$ has generators $\{a,b\}$ and no relations, and the abelianization has the same generators, but now they commute. 

\textbf{Length and exponent-sum.} Every element $g$ in a group $G$ generated by $\Sigma$ has a length $|g|_{\Sigma}$, which is the word length $|.|$ of a shortest word $w$ representing $g$ in $G$. For example, the length of $aba^{-1}$ in $\mbb{Z}^2$ is $|aba^{-1}|_{\Sigma}=|b|=1$.

For any generator $x$ of $F(\Sigma)$, the map $|| . ||_x: F(\Sigma) \to \mbb{Z}$ represents the \emph{exponent-sum} of $x$ in a word $w$; that is, $||w||_x =|w|_x -|w|_{x^{-1}}$, so for example $||xyx^{-1}y^2||_x=0$. One can define the exponent-sum of a generator $x$ in an element $g$ for certain (but not all) groups beyond free groups, and then we use the same notation $||g||_x$.
The length and the image under abelianization are well-defined for any element in any finitely generated group. However, the situation regarding the exponent-sum (of a generator) is more complicated. 
For example, if $H$ is a group and $x\in H$ a generator of order $5$, then one may claim the exponent-sum of $x$ in the element $x^3$ to be $3$; but $x^3=x^{-2}$ in $H$, and in $x^{-2}$ the exponent-sum of $x$ appears to be $-2$.

\textbf{Formal languages and constraints.}
As before, suppose $G$ is a group generated by $S(=\Sigma \cup \Sigma^{-1})$, and let $\pi: S^* \rightarrow G$ be the natural projection from the free monoid $S^*$ generated by $S$ to $G$, taking a word over the generators to the element it represents in the group. A language over $S$ is \emph{regular} if it is recognised by a finite state automaton, as is standard.
%When $w$ is a word over $S$, and we write $|w|$ to denote the length of $w$.

\begin{definition}\label{def:rat_reg_con}
\ \begin{enumerate}
\item[(1)] 
A subset $L$ of $G$ is {\em recognisable} if the full preimage
$\pi^{-1}(L)$ is a regular subset of $S^*$.
\item[(2)] A subset $L$ of $G$ is {\em rational} if $L=\pi(L')$, where $L'$  is a regular subset of $S^*$.
%\item[(3)] A regular subset $L'$ of $S^*$ is
%\emph{quasi-isometrically embedded} (q.i. embedded) in $G$ if there exist
%$\lambda \geq 1$ and $\mu \geq 0$ such that, for any $w \in L'$,
%$|\pi(w)|_G \geq \frac{1}{\lambda}|w|-\mu.$
%\item[(4)] A rational subset $L$ of $G$ is \emph{quasi-isometrically
%embeddable} (q.i. embeddable) in $G$ if there exists a quasi-isometrically embedded
%regular subset  $L'$ of $S^*$ such that $\pi(L')=L.$
\end{enumerate}
\end{definition}
It follows immediately  that recognisable subsets of $G$ are
rational.

\subsection{Word equations with constraints in groups.} Word equations in groups are essentially the same as word equations in monoids, the only difference consisting in the fact that inverses, of both elements and varialbes, are allowed.

 For a group $G$, a \textit{finite system of equations} in $G$ over the variables $\Omega$ is a
		finite subset $\mathcal{E}$ of of the free product $G \ast F(\Omega)$, where
		\(F(\Omega)\) is the free group on \(\Omega\). If
		\(\mathcal{E} = \{w_1, \ \ldots, \ w_n\}\), then a \textit{solution} to the
		system \(w_1 = \cdots = w_n = 1\) is a homomorphism \(\phi \colon G \ast F(\Omega) \to
		G\), such that \(\phi(w_1) = \cdots =\phi (w_n) = 1_G\), where
		 \(\phi(g) = g\) for all \(g \in G\). If $\mathcal{E}$ has a solution, then it is \emph{satisfiable}.
		 
	\begin{example} Consider the system $\mathcal{E}=\{w_1, w_2\}\subset F(a,b) \ast F(X_1, X_2)$ over the free group $F(a,b)$, where $w_1=X_1^2(abab)^{-1}$, $w_2=X_2X_1X_2^{-1}X_1^{-1}$; we set $w_1=w_2=1$, which can be written as $X_1^2=abab, X_2X_1=X_2X_1$. The solutions are $\phi(X_1)=ab, \phi(X_2)=(ab)^k$, $k \in \mbb{Z}$.
	\end{example}

For a group $G$, we say that systems of equations over $G$
are \emph{decidable} over $G$ if there is an algorithm to determine whether any
given system is satisfiable.  The question of decidability of (systems of) equations is often called the
{\em Diophantine Problem} for $G$, and will be denoted here by (\textbf{WordEqnGp$(G)$}).

 The \emph{Diophantine Problem with rational} constraints, denoted (\textbf{WordEqnGp, RAT}), or \emph{recognisable constraints}, denoted (\textbf{WordEqnGp, REC}), asks about the existence of solutions to $\mathcal E$, with some of the variables restricted to taking values in specified rational or recognisable sets, respectively. 
 
 One can attach similar types of non-rational constraints in the group setting as for word equations in free monoids. We illustrate the kind of problems concerning group equations with non-rational constraints by an easy example, where each question represents one of the types of constraints that we are interested in. For each case the problem is about the existence of an algorithm to decide the satisfiability of an equation with constraints.

\begin{example} Consider the equation 
\begin{equation} \label{eqngp}
XaY^2bY^{-1}=1
\end{equation}
 over variables $X,Y$ 
in the free group on two generators $F(a,b)$ with length function $|.|=|.|_{\{a,b\}}$ and abelianization $(\mbb{Z}^2,+)$.
\begin{enumerate}
\item (\textbf{WordEqnGp, LEN}) asks about the existence of solutions whose lengths satisfy some given linear system of equations over the integers.

\medskip

An instance of this problem is: decide whether there are any solutions $(x,y)$ to (\ref{eqngp}) such that $|x|=|y|+2$; the answer is yes, since $(x,y)=(b^{-1}a^{-1}, 1)$ is a solution with $|x|=|y|+2$.
\item  (\textbf{WordEqnGp, EXP-SUM}) asks about the existence of solutions where the exponent-sums of generators satisfy some given linear system of equations over the integers.

\medskip

An instance of this problem is: decide whether there are any solutions $(x,y)$ to (\ref{eqngp}) such that $|x|_a=2|y|_a+|y|_b$ and $|x|_b=3|y|_b$; the answer is `no' by solving a basic linear system over the integers with variables $|x|_a, |x|_b, |y|_a, |y|_b$ and we leave this as an exercise.

\item (\textbf{WordEqnGp, \ab}) asks about the existence of solutions in $G$ whose abelianisation satisfy some given system of equations in $G^{\ab}$.

\medskip

An instance of this problem is: decide whether there are any solutions $(x,y)$ to (\ref{eqngp}) such that $\ab(x)=3\ab(y)$ (we use additive notation for $\mbb{Z}^2$); the answer is no, since $\ab(xax^2by^{-1})=\ab(x)+\ab(y)+(1,1)=(0,0)$ together with $\ab(x)=3\ab(y)$ lead to $4\ab(y)=(-1,-1)$, which is not possible in $\mbb{Z}^2$.
\end{enumerate}
\end{example}

Several other constraints can be imposed, such as requiring that the solutions belong to specified context-free sets, or that they can be compared via a lexicographic order. See \cite{CiobanuEvettsLevine} for details.

Positive results have been established for virtually abelian groups, that is, groups that have abelian subgroups of finite index. 

\begin{theorem}[\cite{CiobanuEvettsLevine}]\label{thm:main}
	In any finitely generated virtually abelian group, it is effectively decidable whether a finite system of equations with the following kinds of constraints has solutions:
	\begin{itemize}
	\item[(i)] linear length constraints (with respect to any weighted word metric),
	\item[(ii)] abelianisation constraints,
	\item[(iii)] context-free constraints,
	\item[(iv)] lexicographic order constraints.
\end{itemize}
	\end{theorem}

The main undecidability results for groups have so far been obtained with respect to (\textbf{WordEqnGp, \ab}):

% which we also call \emph{abelianization constraints}, and denote by $\ab$; these require that the `abelian form' of the solutions, where any two symbols commute, satisfies a given set of equations as well (equivalently, the constraints can be viewed as equations in the abelianization of $G$). We denote this problem by $\mc{DP}(G,\ab)$, and note that this can be seen as introducing an abelian predicate or relation to the existential theory of $G$. %\marginpar{\ag{I may delete this sentence, since it is arguable that the size of the abelianization is such a factor: for example a free abelian group can have very large abelianization rank. Also, at this point removing it leaves us at exactly 15 pages} It turns out that the size of the abelianization of the group $G$ is a major factor: when the abelianization is infinite and `sufficiently large', as for (non-abelian) partially commutative groups, $\mc{DP}(G,\ab)$ is undecidable, while when the abelianization is finite we can reduce abelian constraints to recognisable ones and $\mc{DP}(G,\ab)$ is decidable. } 
 
 \begin{thmletter}[\cite{CiobanuGarreta}]%(Theorem \ref{thm:DP_RAAGS})
	Let $G$ be a partially commutative group (or right-angled Artin group) that is not abelian. Then (\textbf{WordEqnGp, \ab}) is undecidable.
\end{thmletter}

 \begin{thmletter}[\cite{CiobanuGarreta}]
Let $G$ be a hyperbolic group with abelianisation of torsion-free rank $\geq 2$. Then (\textbf{WordEqnGp, \ab}) is undecidable. 
 \end{thmletter}

The key approach for establishing Theorems A and B is the interleaving of algebra and model theory. The main tool used to prove undecidability results is \emph{interpretability} using disjunctions of equations (i.e.\ positive existential formulas), or \emph{PE-interpretability}, which allows one to translate one structure into another and to reduce the Diophantine Problem from one structure to the other. The main reductions are to the Diophantine Problem in the ring of integers, which is a classical undecidable problem (Hilbert's 10th problem). %We give a technical result in Section \ref{sec:technical} which enables us to encode the Diophantine Problem over the integers into the Diophantine Problem with abelianization constraints in a group, assuming the group satisfies certain model-theoretic properties. This technical result (Lemma \ref{l: technical_lemma}) is then applied to partially commutative groups, and in work in progress also to hyperbolic groups with large abelianization or more general graph products. 

As for free monoids, the starting point is B\"uchi and Senger's paper \cite{RichardBuchi1988}, where they show that (positive) integer addition and multiplication can be encoded into word equations in a free semigroup, when requiring that the solutions satisfy an abelian predicate; by the undecidability of Hilbert's 10th problem, such enhanced word equations are also undecidable.
 %the Diophantine problem $\mc{DP}(\Sigma^*, \ab)$ is also undecidable.
  In order to encode the multiplication in the ring $(\mbb{Z}, \oplus, \odot)$ within a (semi)group, one needs two `independent' elements, which can be taken to be two of the free generators if the (semi)group is free. In non-free groups one requires the existence of two elements that play a similar role; however, it is not enough to pick two generators, these elements also need to satisfy additional properties with respect to the abelianization. % project to sufficiently independent  need to also be generators in the infinite part of abelianization of the group, and we therefore name these elements \emph{abelian-primitive};
   Finding such a pair of elements, called \emph{abelian-primitive} in \cite{CiobanuGarreta}, is difficult or impossible for arbitrary groups, but they can be found in (non-abelian) partially commutative groups and hyperbolic groups with `large' abelianisation. 
 
Finally, Theorem \ref{thmC} provides an interesting contrast to the previous results, in that positive outcomes for (\textbf{WordEqnGp, \ab}) hold when $G$ has finite abelianization. In this case the (\textbf{WordEqnGp, \ab}) can be reduced to the (\textbf{WordEqnGp}) with recognisable constraints, and we showcase groups for which these problems are decidable. 

\begin{thmletter}[\cite{CiobanuGarreta}]\label{thmC}%(Theorems \ref{thm:finiteab}, \ref{thm:hyp_finiteab}, \ref{thm:graph_prod})
Let $G$ be a group where the (\textbf{WordEqnGp, REC}) is decidable. Then (\textbf{WordEqnGp, \ab}) is decidable. In particular, this holds if:
\begin{enumerate}
\item[1.] $G$ is a hyperbolic group with finite abelianization.
\item[2.] $G$ be a graph product of finite groups, such as a right-angled Coxeter group.
\end{enumerate}
\end{thmletter}

We finish this summary of results with some illustrative examples.% that are a consequence of Theorems A, B and C.
\begin{example}
\begin{itemize}
\item[1.] The graph product $G$ of finite cyclic groups $G_i$ based on 
 \[
  \xymatrix{\stackrel{G_1}{\bullet}\ar@{-}[r] &\stackrel{G_2}{\bullet} \ar@{-}[r]&\stackrel{G_3}{\bullet}\ar@{-}[r] &\stackrel{G_4}{\bullet} }
 \]
is $\langle x_1,x_2,x_3,x_4 \mid [x_1,x_2]=[x_2,x_3]=[x_3,x_4]=1, x_1^2=x_2^4=x_3^6=x_4^8=1 \rangle,$
has finite abelianisation $G^{\ab}$, and so (\textbf{WordEqnGp, \ab}) is decidable by Theorem C.

%\item $G^{\ab}$ is the direct product $C_2 \times C_4 \times C_6 \times C_8$ of finite cyclic groups of orders $2, 4, 6, 8$, respectively.

\medskip

\item[2.] The surface group $G= \langle a,b,c,d \mid aba^{-1}b^{-1}cdc^{-1}d^{-1}=1\rangle$
is hyperbolic and the abelianisation $G^{\ab}$ of $G$ is $\mathbb{Z}^4$, so (\textbf{WordEqnGp, \ab}) is undecidable by Theorem B.

\medskip

\item[3.] The free product $G= \langle x_1,x_2,x_3,x_4 \mid x_1^2=x_2^4=x_3^6=x_4^8=1 \rangle$ of finite groups 
  is hyperbolic and has finite abelianisation, so (\textbf{WordEqnGp, \ab}) is decidable by Theorem C.

 \end{itemize}
 \end{example}
 % Section \ref{sec:technical} provides the technical tools (Lemma \ref{l: technical_lemma}) required to encode the Diophantine Problem over the integers into a group Diophantine Problem with abelianization constraints, if the group satisfies certain model-theoretic properties. We use this technical section to establish the most involved result of the paper (Theorem \ref{thm:DP_RAAGS}) in Section \ref{sec:RAAGs}, namely, that the $\mc{DP}$ with abelianization constraints is undecidable in (non-abelian) partially commutative groups. In the final section we outline several of the many possible directions of future research that naturally stem from this paper.

  \section{Conclusions}

Word equations in free monoids are an established research direction in theoretical computer science, and imposing constraints on the solutions is of interest for both theoretical and practical reasons. Most results tend to lead to undecidability, and some outstanding open problems remain. In this survey we highlighted a rare positive result, Theorem
\ref{mainGeorg}, for a type of constraints that require none of the solutions to satisfy certain properties.
 
The same type of problems can be posed in the infinite groups for which the satisfiability of equations is known to be decidable, such as virtually abelian, partially commutative groups (the group counterparts of trace monoids) and more generally, graph products of groups, as well as hyperbolic groups.  The motivation to study non-rational constraints is twofold: first, explore extensions of word equations that have not been systematically studied for groups before, and second, develop algebraic and model-theoretic tools that can complement the combinatorial techniques used for solving word equations. As for monoids, most results so far show undecidability, but for virtually abelian groups and certain groups with finite abelianisation several problems turn out to be decidable. Many interesting questions are yet to be explored.

%ultimately augment our understanding of the Diophantine Problem in both groups and semigroups.

The systematic crossover of non-rational constraints from computer science into algebra has begun only recently, and there are undoubtedly many more avenues to explore by translating questions on free monoids, where decades of literature exists, into groups. Once in the world of groups, there are numerous tools to tackle these problems which are not always available for free monoids, tools coming from geometry, topology and algebra.
We expect that the progress and ideas from algebra can inform the work in computer science. 
\bibliographystyle{amsalpha}
\def\cprime{$'$} \def\cprime{$'$}

\end{document}